\magnification=\magstep1
\hsize=16.5 true cm 
\vsize=23.6 true cm
\font\bff=cmbx10 scaled \magstep1
\font\bfff=cmbx10 scaled \magstep2
\font\bffg=cmbx10 scaled \magstep3

\parindent0cm
\def\cl{\centerline}           %
\def\bp{\bigskip}              %
\def\mp{\medskip}              %
\def\sp{\smallskip}            %
\font\boldmas=msbm10           %
\def\Bbb#1{\hbox{\boldmas #1}} %
\def\Q{\Bbb Q}                 %
\def\R{\Bbb R}                 %
\def\N{\Bbb N}                 %
\def\Z{\Bbb Z}                 %
\def\P{\Bbb P}                 %
\centerline{\bffg Many non-embeddable infinite groups}
\bigskip
\centerline{\bfff Gerald Kuba}
\bp
{\bff 1. Introduction}
\mp
Write $\,|M|\,$ for the cardinal number (the {\it size}) of a set $\,M\,$.
If $\,G\,$ is a group then $\,|G|\,$ is the {\it order} of $\,G\,$.
Throughout the paper, $\,{\bf b}_\alpha\,$ denotes 
the $\alpha$th beth. 
(So $\,{\bf b}_0=\aleph_0=|\N|\,$ and
$\;{\bf b}_{\alpha}\,=\,\sup\,\{\,2^{{\bf b}_\beta}\;|\;\beta<\alpha\,\}\;$
for every ordinal $\,\alpha>0\,$. 
For example, $\,{\bf b}_1=2^{{\bf b}_0}=|\R|\,$
and $\,{\bf b}_2=2^{{\bf b}_1}\,$ is the size of the power set of $\,\R\,$
and $\,{\bf b}_3=2^{{\bf b}_2}\,$ and so on.)
In particular, $\,{\bf b}_{\omega_1}\,$ is 
the smallest beth greater than uncountably many beths.
Equivalently, $\,{\bf b}_{\omega_1}\,$ is the 
least strong limit greater than uncountably many strong limits.
(A cardinal number $\,\lambda>0\,$ is a {\it strong limit} 
if $\,2^\kappa<\lambda\,$ for every cardinal number 
$\,\kappa<\lambda\,$.) 
Consider the following two statements for an arbitrary infinite cardinal 
number $\,\kappa\,$.
\sp
(1.1)\quad {\it There exist $\,2^\kappa\,$ pairwise non-isomorphic 
groups of order $\,\kappa\,$.} 
\sp
(1.2)\quad {\it There exist $\,2^\kappa\,$ pairwise non-embeddable
groups of order $\,\kappa\,$.} 
\sp
Trivially, (1.2) implies (1.1).
It goes without saying that in 
(1.1) and hence in (1.2) the number $\,2^\kappa\,$ cannot be exceeded.
As proved by Kurosh [2], (1.1) is true 
for every infinite cardinal $\,\kappa\,$.
Shelah proved in 1974 (see [5] Theorem 2.1) that 
(1.2) is true for every infinite cardinal $\,\kappa\,$
even in the realm of {\it torsion-free abelian groups}. 
As proved by Macintyre and Shelah in 1976 (see [3] Theorem 8), 
(1.2) is true for every {\it regular} cardinal $\,\kappa>\aleph_0\,$
in the realm of {\it universal locally finite groups}.
Since all finite groups embed into each universal locally finite group,
our first goal is to write down a short proof of the following theorem. 
For the sake of simplicity, a group is {\it F-universal} 
if all finite groups embed into it.
\mp
{\bf Theorem 1.} {\it For every infinite cardinal $\,\kappa\,$ 
there exist $\,2^\kappa\,$ 
pairwise non-embeddable 

F-universal groups of order $\,\kappa\,$.}
\mp
While in Theorem 1 non-embeddability is guaranteed only for 
equipotent groups, to make sure 
that small groups are never embeddable in large groups
is a challenge in proving the following theorem
which generalizes Theorem 1. 
\mp
{\bf Theorem 2.} {\it If $\,K\,$ is a set of infinite cardinal numbers
with $\,|K|={\bf b}_{\omega_1}\,$ then a family $\,{\cal F}\,$ of 
pairwise non-embeddable F-universal groups exists such that 
$\,{\cal F}\,$ contains $\,2^\kappa\,$ groups of order $\,\kappa\,$ 
for every $\,\kappa\in K\,$.}
\mp
The proofs of Theorems 1 and 2 are not self-contained because 
we apply [3] Theorem 8, a deep theorem proved
by methods from model theory and first-order logic.
For the proof of Theorem 2 we also need
a specific construction which proves the following theorem.
Since this construction is carried out in a vivid and 
elementary way, we can provide a self-contained proof of 
(1.2) in the realm of torsion-free groups
for the range $\,\aleph_0\leq\kappa\leq{\bf b}_{\omega_1}\,$.
\mp
{\bf Theorem 3.} {\it There exists a family $\,{\cal H}\,$ of 
pairwise non-embeddable torsion-free groups such that 
$\,{\cal H}\,$ contains $\,2^\kappa\,$ groups of order $\,\kappa\,$ 
for every infinite cardinal $\,\kappa\leq{\bf b}_{\omega_1}\,$.} 
\mp
{\it Remark.} While Theorem 3 covers 
precisely $\,\aleph_1\,$ beths, 
it is possible that it is a statement 
about more than $\,\aleph_1\,$ alephs. 
As Theorem 2, 
Theorem 3 could even be a statement about $\,{\bf b}_{\omega_1}\,$ cardinals.
In fact, 
{\it it is consistent with {\rm ZFC} set theory that $\,{\bf b}_\alpha\,$
cardinal numbers lie 
below $\,{\bf b}_\alpha\,$ for every limit ordinal $\alpha$.}
A challenge in proving Theorem 3 is also that 
{\it it is consistent with {\rm ZFC} that there exist 
$\,{\bf b}_{n+1}\,$ cardinals $\,\kappa\,$ with 
$\,{\bf b}_{n}<\kappa<{\bf b}_{n+1}\,$ for every even integer $\,n\geq 0\,$.}
(An explanation for these consistency results is given in Section~7.)
\eject
\bp
{\bff 2. Preparation of the proofs}
\mp\sp
For the sake of simplicity, groups are {\it incomparable} if 
they are {\it pairwise non-embeddable}. 
Let $\,{\cal G}\,$ be a set of nontrivial groups
with $\,|{\cal G}|>1\,$.
Then $\,\prod^*\!{\cal G}\,=\,\prod^*_{G\in {\cal G}}G\,$
denotes the {\it free product} 
of the groups in $\,{\cal G}\,$ defined in the usual way (see [4]).
(Since a free product is invariant under replacing 
a factor with an isomorphic one, we need not consider
free products $\,\prod^*_{i\in I}G_i\,$
where the mapping $\,i\mapsto G_i\,$ is not injective.)
The group $\,\prod^*\!{\cal G}\,$ 
is {\it infinite} and non-abelian.  More precisely, 
\mp\sp
(2.1)\quad $\;|\prod^*\!{\cal G}|\,=\,\max\big\{\,\aleph_0,\,|{\cal G}|,\,
\sup\{\,|G|\;|\;G\in{\cal G}\,\}\,\big\}\;.$ 
\mp\sp
In particular, $\;|\prod^*\!{\cal G}|\,=\,|{\cal G}|\;$ 
if $\;\aleph_0\leq|G|\leq|{\cal G}|\;$ for every $\,G\in{\cal G}\,$.
In case that $\,{\cal G}=\{G,H\}\,$ (and $\,G\not= H\,$)
we simply write $\,G*H\,$ instead of $\,\prod^*\{G,H\}\,$.
Thus $\;|G*H|=\max\{|G|,|H|\}\,$ if $\,G\,$ or $\,H\,$ is infinite.
\mp\sp
It is well-known (see [4] Ex.~11.65) that a nontrivial {\it direct product}
of groups cannot be a free product. For our purpose 
it is enough to apply the following lemma.
\mp
{\bf Lemma 1.} {\it Let $\,G,H\,$ be groups with $\,|G|,|H|>1\,$.  
If $\,G\,$ is abelian then the direct product $\,G\times H\,$
of the groups $\,G\,$ and $\,H\,$ is not a free product of groups.}
\mp
{\it Proof.} Naturally, the center 
of $\,G\times H\,$ contains the abelian factor $\,G\,$ 
and $\,|G|>1\,$. On the other hand, 
free products (with at least two factors) are always centerless, {\it q.e.d.}
\mp\sp
In the following let us call a group $\,G\,$ {\it freely indecomposable}
or, more simply, a {\it FI-group} if $\,|G|>1\,$ and 
$\,G\,$ is not isomorphic with the free group $\,\Z\,$
or with a free product of groups. 
In particular, 
torsion groups and abelian groups not isomorphic with $\,\Z\,$ 
are FI-groups. (A free group of rank $\,\kappa>1\,$ 
is isomorphic to a free product of $\,\kappa\,$ distinct copies of $\,\Z\,$
and hence no free group is a FI-group.)
By Lemma 1, if $\,G\,$ is an abelian group 
and $\,|G|>1\,$ then $\,G\times H\,$ 
is an FI-group for every group $\,H\,$ with $\,|H|>1\,$.
The following lemma is an immediate consequence of the famous Kurosh theorem 
[4] 11.55 about subgroups of free products.
\mp
{\bf Lemma 2.} {\it If $\,G\,$ is a FI-group and 
$\,{\cal H}\,$ is a family of FI-groups with $\,|{\cal H}|>1\,$ 
and there exists a monomorphism from 
$\,G\,$ into $\,\prod^*\!{\cal H}\,$ then 
there exists a monomorphism from 
$\,G\,$ into some group $\,H\in {\cal H}\,$.}
\mp
The following immediate consequence of Lemma 2 is 
essential for our proofs.
\mp
{\bf Corollary 1.} {\it Let $\,{\cal G}\,$ and $\,{\cal H}\,$ 
be families of FI-groups with $\,|{\cal G}|>1\,$ and $\,|{\cal H}|>1\,$. 
If some group in $\,{\cal G}\,$ cannot be embedded in any group 
in $\,{\cal H}\,$ and some group in $\,{\cal H}\,$ cannot be embedded
in any group in $\,{\cal G}\,$  then the groups 
$\,\prod^*\!{\cal G}\,$ and $\,\prod^*\!{\cal H}\,$ are 
incomparable.} 
\mp\sp
For the proofs of Theorems 1 and 3 
we also need an operator $\,{\cal Y}\,$ 
on the class of all infinite sets.
\mp
(2.2)\quad {\it If $\,S\,$ is an infinite set then 
$\,{\cal Y}(S)\,$ is a family  
of subsets of $\,S\,$ such that 
$\;|X|=|S|\;$ for every $\,X\in{\cal Y}(S)\,$ and 
$\,|{\cal Y}(S)|=2^{|S|}\,$ and $\;X\not\subset Y\;$ 
whenever $\,X,Y\in{\cal Y}(S)\,$ are distinct.}
\mp
For example, write $\,S=S_1\cup S_2\,$ 
with $\,|S_1|=|S_2|\,$ and $\,S_1\cap S_2=\emptyset\,$.
Let $\,f\,$ be a bijection from $\,S_1\,$ onto 
$\,S_2\,$ and put $\;{\cal Y}(S)\,=\,
\{\,T\cup(S_2\setminus f(T))\;|\;T\subset S_1\,\}\,$.
With this definition, (2.2) is obviously true. 
\vfill\eject
\bp
{\bff 3. Countable groups}
\mp\sp
In order to settle the special case $\,\kappa=\aleph_0\,$
in Theorem 1 we consider specific subgroups of $\,\Q\,$.
Let $\,\P\,$ be the set of all primes and for $\,P\subset\P\cup\{1\}\,$ 
put 
$$\Delta[P]\,=\,
\Big\{\,{m\over p_1\cdots p_n}\;\Big|\;m\in\Z\;\land\;
p_1,...,p_n\in P\,\Big\}$$                                  
In particular, $\,\Delta[\{1\}]=\Z\,$ and $\,\Delta[\P]=\Q\,$.
Instead of the additive group $\,\Delta[P]\,$ 
we also consider the isomorphic multiplicative group 
$\,\Gamma[P]\,:=\,\{\,e^x\;|\;x\in\Delta[P]\,\}\,$ and put
\sp
\cl{$\;\Gamma_p\,:=\,\Gamma[\{p\}]\,=\,
\{\,e^{np^k}\;|\;n,k\in\Z\,\}\;$\quad for every $\,p\in\P\cup\{1\}\,$.}  
\mp
The {\it additively} written abelian group $\,\Delta[P]\,$ has 
an essential property 
that motivates the following definition 
for arbitrary groups in {\it multiplicative} notation.
\mp
(3.1)\quad {\it If $\,G\,$ is a group and $\,p\,$ 
is a prime number and $\,a\in G\,$
then $\,a\,$ {\it is $p$-divisible in $\,G\,$}
if and only if for each $\,n\in\N\,$ the equation 
$\,x^{p^n}=a\,$ has a solution $\,x\in G\,$.}
\mp
So if $\,P\subset\P\,$ and $\,q\in\P\setminus P\,$ 
then every element of $\,\Gamma[P]\,$ is $p$-divisible in $\,\Gamma[P]\,$
for every $\,p\in P\,$ and 
no element of $\,\Gamma[P]\setminus\{1\}\,$ 
is $q$-divisible in $\,\Gamma[P]\,$.
\mp\sp
Let $\,F\,$ be the group of all permutations $\,f\,$
of $\,\N\,$ where $\,\{\,n\in\N\;|\;f(n)\not=n\,\}\,$ 
is a finite set. Obviously, 
the group $\,F\,$ is F-universal by Cayley's theorem.
Of course, $\,F\,$ is a countably infinite torsion group.
The first step in proving Theorem 1 is to verify the following statement
using (2.2).
\mp
(3.2)\quad {\it The $\,2^{\aleph_0}\,$ F-universal 
countable FI-groups 
\sp
\cl{$\Gamma[P]\times F\;\;\;\;
(P\in{\cal Y}(\P))$}
\sp
are incomparable.} 
\mp
{\it Proof.} Let $\,P,P'\subset\P\,$ 
with $\,P\not\subset P'\,$ and assume that the group
$\,\Gamma[P]\times F\,$ is embeddable in the group 
$\,\Gamma[P']\times F\,$.
Then there is a monomorphism $\,\varphi\,$ from $\,\Gamma[P]\,$ 
to $\,\Gamma[P']\times F\,$.
Choose $\,p\in P\setminus P'\,$ and $\,a\in\Gamma[P]\setminus\{1\}\,$.
Then $\,a\,$ is $p$-divisible in the group $\,\Gamma[P]\,$
and of infinite order.
But $\,(b,f)=\varphi(a)\,$ is either of finite order 
(in case $\,b=1\,$)
or (in case $\,b\not=1\,$)
not $p$-divisible in the group $\,\Gamma[P']\times F\,$ 
and hence not $p$-divisible in the group 
$\,\varphi(\Gamma[P])\,$. 
This contradiction finishes the proof of (3.2).
\mp
For the proof of Theorems 1 and 2 we also need a lemma 
about $p$-divisibility.
For $\,P\subset\P\,$ let $\,{\bf G}[P]\,$ be the class of all 
groups $\,G\,$ such that if $\,a\in G\,$ then 
$\;a=a_1a_2\cdots a_n\;$ where each $\,a_i\in G\,$
is either of finite order or $p$-divisible in $\,G\,$ 
for some $\,p\in P\,$.
Obviously, for every 
$\,P\subset\P\,$ the class $\,{\bf G}[P]\,$ is closed under free products.
\mp
{\bf Lemma 3.} {\it If $\,P\subset\P\,$ and 
$\,Q\subset(\P\cup\{1\})\setminus P\,$ and  
$\,G,H\,$ are groups with 
$\,G\in{\bf G}[P]\,$ and 
$\,G\,$ is embeddable in  
$\,\Gamma[Q]\!\times\!H\,$ then $\,G\,$ is embeddable in $\,H\,$.}
\mp
{\it Proof.} Let $\,\varphi\,$ be a monomorphism from 
$\,G\,$ to $\,\Gamma[Q]\!\times\!H\,$
and $\,a\in G\,$. 
Then $\;a=a_1a_2\cdots a_n\;$ where each $\,a_i\in G\,$
is either of finite order or 
$p$-divisible in $\,G\,$ for some $\,p\in P\,$.
Put $\,\varphi(a_i)=(b_i,c_i)\,$. 
Then $\,\varphi(a)=(b_1,c_1)\cdots(b_n,c_n)=(b_1\cdots b_n,c_1\cdots c_n)\,$
and every $\,(b_i,c_i)\,$ is either of finite order or 
$p$-divisible in $\,\varphi(G)\,$ and hence in $\,\Gamma[Q]\!\times\! H\,$
for some $\,p\in P\,$.
Since $\,Q\cap P=\emptyset\,$ this is only possible
if $\,b_i=1\,$ for every $\,i\,$. 
Hence $\,\varphi(a)=(1,c)\,$ with $\,c\in H\,$ and 
therefore $\,\varphi(G)\subset \{1\}\!\times\! H\,$.
This is enough since the groups 
$\,\{1\}\!\times\! H\,$ and $\,H\,$ are isomorphic, {\it q.e.d.}
\vfill\eject

\bp
{\bff 4. Proof of Theorem 1}
\mp\sp
There are two ways to settle Theorem 1. 
The first one is very short and 
applies Shelah's theorem about abelian groups from [5]. 
The second one applies [3] Theorem 8 and
verifies also an additional condition we need for the proof of Theorem 2. 
\mp
The short proof of Theorem 1 reads as follows. 
By [5] Theorem 2.1, for every infinite cardinal $\,\kappa\,$ there 
exists a family $\,{\cal A}_\kappa\,$ 
of $\,2^{\kappa}\,$ incomparable torsion-free abelian groups of
order $\,\kappa\,$. Consider the countable permutation group 
$\,F\,$ from Section 3. Similarly as in the proof of (3.2) 
the $\,2^{\kappa}\,$ F-universal groups 
$\,A\times F\;(A\in{\cal A}_\kappa)\,$ obviously settle Theorem 1. 
\mp
Considering (3.2), the second proof is the proof of the following statement.
\mp
(4.1)\quad {\it For every uncountable cardinal number $\,\kappa\,$  
there exist $\,2^\kappa\,$ pairwise non-embeddable 
F-universal groups $\,G\,$ of order $\,\kappa\,$
such that if $\,a\in G\,$ is of infinite order then 
$\,a\,$ is not $p$-divisible for any prime $\,p>3\,$.}
\mp
{\it Proof.} As usual, if $\,\kappa\,$ is an infinite cardinal then 
$\,\kappa^+\,$ denotes the smallest cardinal 
greater than $\,\kappa\,$. (Notice that $\,\kappa^+\,$
is always a {\it regular} cardinal.)
By [3] Theorem 8, for every regular $\,\kappa>\aleph_0\,$ 
there exists a family $\,{\cal F}_\kappa\,$
of $\,2^{\kappa}\,$ incomparable FI-groups of
order $\,\kappa\,$ which are F-universal torsion groups. 
So let $\,\kappa\,$ be singular.
Then for some regular, infinite $\,\mu<\kappa\,$ 
we have $\;\kappa\,=\,\sup\Lambda\;$
for some set $\,\Lambda\,$ of size $\,\mu\,$ where each 
$\,\lambda\in\Lambda\,$ is a {\it regular} cardinal
with $\,\mu^+<\lambda<\kappa\,$. 
\sp
Since $\,\mu^+>\aleph_0\,$ and 
$\,|\Lambda|=\mu<2^\mu\leq 2^{(\mu^+)}=|{\cal F}_{\mu^+}|\,$,
we can define an injective mapping $\,\lambda\mapsto G_\lambda\,$
from $\,\Lambda\,$ to $\,{\cal F}_{\mu+}\,$.
Then $\;G_\lambda\,(\lambda\in\Lambda)\;$
are incomparable FI-groups of order $\,\mu^+\,$.
Since $\,|{\cal F}_\lambda|=2^\lambda\,$ for each $\,\lambda\in\Lambda\,$, 
we can choose a set $\,{\cal F}'_\lambda\subset{\cal F}_\lambda\,$
of size $\,\lambda\,$ for each $\,\lambda\in\Lambda\,$.
Trivially, the families $\,{\cal F}'_\lambda\,(\lambda\in \Lambda)\,$
are pairwise disjoint. 
Now define 
\mp
\cl{${\cal D}\;=\;\{\,\Gamma_2\times((\Gamma_3\!\times\!G_\lambda)*H)\;\,
|\,\;\lambda\in\Lambda\;\land\;
H\in{\cal F}'_\lambda\,\}\;.$}
\mp
In view of Lemma 1, each group in the family $\,{\cal D}\,$ is a FI-group.
The fourth factors ensure that each group in $\,{\cal D}\,$ 
is F-universal.
Since $\;|\Gamma_2\times((\Gamma_3\!\times\!G_\lambda)*H)|
=\max\{\aleph_0,\mu^+,\lambda\}=\lambda\;$
for each $\,\lambda\in\Lambda\,$, the order of 
a group in $\,{\cal D}\,$ is always smaller than $\,\kappa\,$.
Furthermore, $\,|{\cal D}|=\kappa\,$ 
since $\;|\bigcup_{\lambda\in\Lambda}{\cal F}'_\lambda|\,=\,
\sum_{\lambda\in\Lambda}|{\cal F}'_\lambda|\,=\,
\sup\,\{\,\lambda\;|\;\lambda\in\Lambda\,\}\,=\,\kappa\,$.
We claim that the groups in $\,{\cal D}\,$ are incomparable. 
\sp
To prove this let 
$\,\Gamma_2\times((\Gamma_3\!\times\!G_1)*H_1)\,$ and
$\,\Gamma_2\times((\Gamma_3\!\times\!G_2)*H_2)\,$
be two groups in $\,{\cal D}\,$
and assume that 
the first is embeddable in the second. 
The second group has a subgroup isomorphic with 
$\,(\Gamma_3\!\times\!G_1)*H_1\,$ 
and hence $\,(\Gamma_3\!\times\!G_1)*H_1\,$ 
is embeddable in $\,(\Gamma_3\!\times\!G_2)*H_2\,$ by Lemma 3
with $\,P=\{3\}\,$ and $\,Q=\{2\}\,$.
\sp
Since $\,|H_1|>|\Gamma_3\!\times\!G_2|\,$, the FI-group 
$\,H_1\,$ cannot be embedded into the FI-group 
$\,\Gamma_3\!\times\!G_2\,$ and hence $\,H_1\,$ is embeddable 
in $\,H_2\,$ by Lemma 2.
So we obtain $\,H_1=H_2\,$ in case that $\,G_1=G_2\,$.
Therefore the claim is verified by showing that $\,G_1=G_2\,$. 
Since the FI-group $\,\Gamma_3\!\times\!G_1\,$ is not a torsion group, 
it cannot be embedded into the FI-group $\,H_2\,$ and hence 
it must be embeddable in the FI-group 
$\,\Gamma_3\!\times\!G_2\,$ . A fortiori, $\,G_1\,$ is embeddable 
in $\,\Gamma_3\!\times\!G_2\,$
and hence $\,G_1\,$ is embeddable in $\,G_2\,$ by Lemma 3
with $\,P=\emptyset\,$ and $\,Q=\{3\}\,$, 
whence $\,G_1=G_2\,$.
\sp
Finally, using (2.2), 
\mp
\cl{${\cal A}\,=\,\{\,\prod^*\!X\;\,|\,\;X\in{\cal Y}({\cal D})\,\}$} 
\mp
is a family of $\,2^\kappa\,$ 
incomparable groups of order $\,\kappa\,$
in view of (2.1) and Corollary 1. This concludes the proof of (4.1) 
in view of the following lemma which is also essential 
for the proof of Theorem 3. 
\mp
{\bf Lemma 4.} {\it If   
$\,{\cal F}\not=\emptyset\,$ is a family of groups and $\,p\in \P\,$ 
then $\,a\in \prod^*\!{\cal F}\,$ 
is $p$-divisible in $\,\prod^*\!{\cal F}\,$ only if 
$\,a\,$ is conjugate to an element $\,b\,$ of some group
$\,G\in{\cal F}\,$ such that $\,b\,$ is $p$-divisible in $\,G\,$.}
\mp
{\it Proof.} If $\,a\in \prod^*\!{\cal F}\,$ 
then let us call $\,a\,$ {\it good} if there is a group $\,G\in{\cal F}\,$
and $\,x\in G\,$ and $\,y\in\prod^*\!{\cal F}\,$
such that $\,a=yxy^{-1}\,$.
Furthermore, $\,\ell(a)\,$ denotes the length of the word $\,w\,$
obtained from $\,a\,$ by maximal reduction 
via the usual cancellations. (In particular, 
$\,\ell(1)=0\,$.)
We claim that if $\,a\in \prod^*\!{\cal F}\,$ 
then either $\,a\,$ is good or $\,\ell(a^k)\geq k\,$
for every $\,k\geq 2\,$. 
\sp
To verify the claim,  
write $\;a\,=\,a_1\cdots a_nzb_n\cdots b_1\;$
with $\,\ell(a_i)=\ell(b_i)=1\,$ 
and $\,\ell(z)\leq 1\,$ and $\,\ell(a)=2n+\ell(z)\,$.
(So for every group $\,G\in{\cal F}\,$ we always have 
$\,\{a_i,a_{i+1}\}\not\subset G\,$ and  
$\,\{b_i,b_{i+1}\}\not\subset G\,$
and in case that $\,z=1\,$ also 
$\,\{a_n,b_n\}\not\subset G\,$.)
If $\,n=0\,$ then $\,a=z\,$ and hence $\,a\,$ is good. 
Assuming $\,n\not=0\,$ put 
$\;M\,=\,\{\,m\in\N\setminus\{0\}\;|\;
\forall i\leq m\,:\;b_i^{-1}=a_i\,\}\,$. 
If $\,M=\emptyset\,$ then $\,b_1a_1\not=1\,$ or, equivalently,
$\,1\leq\ell(b_1a_1)\leq 2\,$ and hence 
$\;\ell(a^k)=k\ell(a)-(2-\ell(b_1a_1))(k-1)\geq k(2n)-k\geq k\;$ 
for all $\,k\geq 2\,$.
So assume $\,M\not=\emptyset\,$ and put $\,\mu=\max M\,$.
If $\,\mu=n\,$ then $\;a=(a_1\cdots a_n)z(a_1\cdots a_n)^{-1}\;$
is good. If $\,\mu<n\,$ then 
$\,1\leq \ell(b_{\mu+1}a_{\mu+1})\leq 2\,$
and $\;\ell(a^{k})=k\ell(a)-2(k-1)\mu-(k-1)(2-\ell(b_{\mu+1}a_{\mu+1}))
\geq 2k(n-\mu)+2\mu-k\geq k\;$ for each $\,k\geq 2\,$.
This finishes the proof of the claim.
\sp
Now let $\,a\,$ be $p$-divisible in $\,\prod^*\!{\cal F}\,$.
Then for every $\,n\in\N\,$ we can find $\,x_n\in\prod^*\!{\cal F}\,$
such that $\,x_n^{p^n}=a\,$. If $\,x_n\,$ is good 
for some $\,n\in\N\,$ then $\,a\,$ is good 
(since $\,(yxy^{-1})^{p^n}=yx^{p^n}\!y^{-1}\,$) and we are done.
Actually, $\,x_n\,$ must be good for some $\,n\in\N\,$
because otherwise $\,\ell(x_n^{k})\geq k\,$ for every $\,k\geq 2\,$ 
and every $\,n\geq 1\,$ and hence   
we would derive the absurd statement that 
$\,\ell(a)=\ell(x_n^{p^n})\geq p^n\,$ for every $\,n\in\N\,$, 
{\it q.e.d.}
\mp\sp
{\it Remark.} In the proof of (4.1)  
all groups of regular order are FI-groups but 
no group of singular order is FI. 
The proof of Theorem 2 can be carried out much more  easily
if the property FI can be included in (4.1). 
In view of Lemma 3 with $\,P=\{2,3\}\,$ and $\,Q=\{1\}\,$
this can obviously be accomplished by 
replacing the family $\,{\cal A}\,$ 
with the family $\,\{\,\Gamma_1\!\times\!A\;|\;A\in{\cal A}\,\}\,$.        
So (4.1) remains true in the realm of all FI-groups 
in the class $\;{\bf G}[\{1,2,3\}]\,$.
\bp
{\bff 5. Proof of Theorem 3}
\mp
In order to construct appropriate torsion-free groups 
we use the countable abelian groups $\,\Gamma_p\,$ 
from Section 3 as basic building blocks.
\mp
If $\,G\,$ is a group and 
$\,P\,$ is a nonempty set of primes then we 
write $\,\Omega(G)=P\,$ if and only if (i) $\,G\,$ is torsion-free,
(ii) for each $\,p\in P\,$ 
some $\,a\in G\setminus\{1\}\,$ is $p$-divisible in $\,G\,$,
(iii) for each $\,p\in \P\setminus P\,$ no 
$\,a\in G\setminus\{1\}\,$ is $p$-divisible in $\,G\,$, 
(iv) every $\,a\in G\,$ has a representation 
$\,a=a_1\cdots a_m\,$ such that 
for every $\,i\in\{1,...,m\}\,$ there is a prime $\,p_i\in P\,$
with $\,a_i\,$ being $p_i$-divisible in $\,G\,$.
If for a group $\,G\,$ a nonempty set $\,P\subset\P\,$ 
with $\,\Omega(G)=P\,$ does not exist then we write 
$\,\Omega(G)=\emptyset\,$. 
So $\,\Omega(G)\,$ is a well-defined subset of $\,\P\,$
for every group $\,G\,$. 
For example,  $\,\Omega(G)=\emptyset\,$ 
if $\,G\,$ is a free group or a free abelian group or a torsion group.
A pivotal example (with $\,\prod^*\!\{G\}=G\,$) 
and an essential observation read as follows.
\mp
(5.1)\quad {\it If  $\,\emptyset\not=P\subset\P\,$
then $\,\Omega(\prod^*\!\{\,\Gamma_p\;|\;p\in P\,\})=P\,$.}
\mp
(5.2)\quad {\it If  $\,{\cal G}\not=\emptyset\,$ 
is a family of groups $\,G\,$ with $\,\Omega(G)\not=\emptyset\,$
then $\,\Omega(\prod^*\!{\cal G})\,=\,\bigcup_{G\in{\cal G}}\Omega(G)\,$.}
\eject
\mp\sp
(5.1) and (5.2) are true in view of Lemma 4.
The following lemma is an immediate consequence of Lemma 3.
\mp
{\bf Lemma 5.} {\it If $\,p\in\P\,$ 
and $\,G,H\,$ are incomparable groups such that 
$\,\Omega(G),\Omega(H)\,$ are not empty and 
$\,p\not\in\Omega(G)\cup\Omega(H)\,$ 
then the groups $\,\Gamma_p\times G\,$ and
$\,\Gamma_p\times H\,$ are incomparable.}
\mp\sp
If $\,G\,$ is a group with $\,\Omega(G)\not=\emptyset\,$
then $\,\Omega(G)\,$ is the set of all $\,p\in\P\,$
such that for some $\,a\in G\setminus\{1\}\,$ the equation $\,x^{p^n}=a\,$
has a solution in $\,G\,$ for every $\,n\in\N\,$. 
As a trivial consequence, if $\,G\,$ is a subgroup 
of $\,H\,$ and $\,\Omega(H)\not=\emptyset\,$ 
then $\,\Omega(G)\subset\Omega(H)\,$. 
Therefore, the following statement is evident.
\mp
(5.3)\quad {\it If $\,G,H\,$ are groups 
with $\,\Omega(G)\not\subset\Omega(H)\not=\emptyset\,$ 
then $\,G\,$ is not embeddable in $\,H\,$.}
\mp
Of course, the following statement is evident as well.
\mp
(5.4)\quad {\it If $\,H\,$ is a group with $\,\Omega(H)\not=\emptyset\,$ 
and $\,p\in\P\setminus\Omega(H)\,$
then $\,\Omega(\Gamma_p\!\times\!H)=\{p\}\cup\Omega(H)\,$.}
\bp
Now we are ready to prove Theorem 3.  
For abbreviation, a {\it TFFI-group} is a {\it torsion-free FI-group}.
Let $\,\P_1,\P_2,\P_3\,$ denote
three mutually disjoint infinite subsets of $\,\P\setminus\{2,3,5,7\}\,$.
For any ordinal number $\,\alpha\,$ let $\,W_\alpha\,$
denote the set of all ordinals up to $\,\alpha\,$.
(In von Neumann's terminology, $\,W_\alpha=\alpha\cup\{\alpha\}\,$.) 
For the moment we fix an arbitrarily large countable ordinal 
number $\,\eta>0\,$. 
Let $\,\pi\,$ be an injective mapping from 
the set $\;W_\eta\!\times\! W_\eta\,=\,W_\eta^2\;$ to $\,\P_1\,$.
The first step in proving Theorem 3 is to verify 
the following statement.
\mp
(5.5)\quad {\it For every ordinal $\,\xi\leq\eta\,$ a family 
$\,{\cal G}[\xi]\,$ of incomparable TFFI-groups of order 
$\,{\bf b}_\xi\,$ exists
such that $\,|{\cal G}[\xi]|={\bf b}_{\xi+1}\,$
and $\,\emptyset\not=\Omega(G)\subset \P_2\cup\pi(W_{\xi}^2)\,$
for every $\,G\in{\cal G}[\xi]\,$.} 
\mp
In order to prove (5.5) by induction, assume that for an ordinal 
$\,\beta\leq\eta\,$ a family $\,{\cal G}[\alpha]\,$
does the job for every $\,\alpha<\beta\,$. We are done by 
defining a family $\,{\cal G}[\beta]\,$ that does the job.
\mp
If $\,\beta=0\,$ then the assumption is true for all $\,\alpha<\beta\,$
vacuously and hence we have to define $\,{\cal G}[\beta]\,$
appropriately.
This is easily done by considering 
the torsion-free abelian groups $\,\Gamma_p\,$ with $\,p\in\P_2\,$. 
By (5.1) and (5.2) and (5.3), 
the $\,2^{\aleph_0}\,$ groups 
$\;\prod^*\!\{\,\Gamma_p\;|\;p\in P\,\}\;$  
with $\,P\,$ running through $\,{\cal Y}(\P_2)\,$
are incomparable. (Recall the definition (2.2) of $\,{\cal Y}(\cdot)\,$.)
All groups are countably infinite 
in view of (2.1).
While these groups are all torsion-free, 
they are far from being freely indecomposable.
Therefore we define
\mp
\cl{$\;{\cal G}[0]\;=\;
\{\,\Gamma_{\pi(0,0)}\!\times\!{\prod\limits_{p\in P}}\!\!^*\Gamma_p
\,\;|\;\,P\in{\cal Y}(\P_2)\,\}\;.$}
\mp 
In view of (5.1) and Lemma 5, the family $\,{\cal G}[0]\,$ 
consists of $\,{\bf b}_1\,$ incomparable TFFI-groups 
of order $\,{\bf b}_0\,$. By (5.1) and (5.4), if $\,G\in{\cal G}[0]\,$ then 
$\,\emptyset\not=\Omega(G)\subset \P_2\cup\{\pi(0,0)\}\,$.
\mp
Now we assume $\,\beta>0\,$. 
We do not distinguish the cases whether $\,\beta\,$ is a 
successor ordinal or a limit ordinal. 
In both cases we certainly have 
$\;|\bigcup_{\alpha<\beta}{\cal G}[\alpha]|={\bf b}_\beta\,$.
We cannot expect that 
the groups in $\,\bigcup_{\alpha<\beta}{\cal G}[\alpha]\,$
are incomparable. 
However, with 
\mp
\cl{${\cal G}'[\alpha]\;:=\; 
\{\,\Gamma_{\pi(\alpha,\beta)}\!\times\! G\;\,|\,\;
G\in{\cal G}[\alpha]\,\}\;$} 
\mp
the groups in $\;{\cal U}_\beta\,:=\,\bigcup_{\alpha<\beta}{\cal G}'[\alpha]\;$
are incomparable in view of (5.3) and (5.4)
since the primes $\,\pi(\alpha,\beta)\;(\alpha<\beta)\,$
are distinct and do not lie in $\,\bigcup_{\alpha<\beta}\pi(W_\alpha^2)\,$.
We have $\;|{\cal U}_\beta|={\bf b}_\beta\;$ 
and the order of every group in $\,{\cal U}_\beta\,$ 
is less than $\,{\bf b}_\beta\,$.
Finally we define 
\mp
\cl{${\cal G}[\beta]\;=\; 
\{\,\Gamma_{\pi(\beta,\beta)}\times \prod^*\!F\;\,|\,\;
F\in{\cal Y}({\cal U}_\beta)\,\}\,$.}
\mp
This definition concludes the proof of (5.5)
in view of (2.1), (2.2), (5.2), (5.3), (5.4) since
$\;\pi(\beta,\beta)\not\in
\bigcup_{\alpha<\beta}\pi(W_\alpha\!\times\!W_\beta)\,$. 
\mp
Since the {\it countable} ordinal $\,\eta\,$ is arbitrary,
from (5.5) we conclude the following statement
about {\it uncountably many} ordinal numbers.
\mp
(5.6)\quad {\it For every ordinal $\,\eta<\omega_1\,$ a family 
$\,{\cal G}[\eta]\,$ of incomparable TFFI-groups of order 
$\,{\bf b}_{\eta}\,$ exists
such that $\,|{\cal G}[\eta]|={\bf b}_{\eta+1}\,$
and $\,\emptyset\not=\Omega(G)\subset \P_2\cup\P_1\,$
for every $\,G\in{\cal G}[\eta]\,$.} 
\mp
With the help of (5.6) we show that the following is true.
\mp
(5.7)\quad {\it For every ordinal number 
$\,\alpha<\omega_1\,$ there exists a family 
$\,{\cal H}[\alpha]\,$ of TFFI-groups of order $\,{\bf b}_\alpha\,$ 
such that $\,|{\cal H}[\alpha]|={\bf b}_{\alpha+1}\,$
and $\,\bigcup_{\alpha<\omega_1}{\cal H}[\alpha]\,$
consists of incomparable groups    
$\,H\,$ satisfying 
$\,\emptyset\not=\Omega(H)\subset \P\setminus\{2,3,5,7\}\,$.}
\mp
Since there are precisely $\,\aleph_1\,$ ordinals below $\,\omega_1\,$
and since $\,|{\cal Y}(\P_3)|=2^{\aleph_0}\geq\aleph_1\,$,
we can easily settle (5.7) by modifying the families 
$\,{\cal G}[\eta]\,$ in (5.6) as follows. 
Consider the abelian groups $\,\Gamma[P]\,$ from Section 3. 
Of course, $\,\Omega(\Gamma[P])=P\,$ whenever $\,\emptyset\not=P\subset\P\,$. 
Let $\,\alpha\mapsto P_\alpha\,$ be an injective function from 
the set of all countable ordinals to $\,{\cal Y}(\P_3)\,$
and define for each $\,\alpha<\omega_1\,$
\mp
\cl{${\cal H}[\alpha]\;=\; 
\{\,\Gamma[P_\alpha]\times G\;\,|\,\;
G\in{\cal G}[\alpha]\,\}\,$.} 
\mp           
Since $\,\P_3\cap(\P_1\cup\P_2)=\emptyset\,$, 
in view of Lemma 3
we can be sure that for every $\,\alpha<\omega_1\,$ 
the groups in $\,{\cal H}[\alpha]\,$ are incomparable. 
Therefore, (5.7) is settled in view of Lemma 1 and (5.3). 
However, (5.7) is not good enough since the case $\,\alpha=\omega_1\,$
is not included. While omitting the one cardinal 
$\,\kappa={\bf b}_{\omega_1}\,$
in Theorem 3 would not be a tragedy,
Theorem~2 could not be proved if $\,\kappa={\bf b}_{\omega_1}\,$
were not included in Theorem 3. 
Actually, we will need the following modification of (5.7)
for the proof of Theorem 2. 
\mp
(5.8)\quad {\it For every ordinal number 
$\,\alpha\leq \omega_1\,$ there exists a family 
$\,{\cal H}'[\alpha]\,$ of TFFI-groups of order $\,{\bf b}_\alpha\,$ 
such that $\,|{\cal H}'[\alpha]|={\bf b}_{\alpha+1}\,$
and $\,\bigcup_{\alpha\leq\omega_1}{\cal H}'[\alpha]\,$
consists of incomparable groups    
$\,H\,$ satisfying 
$\,\emptyset\not=\Omega(H)\subset \P\setminus\{2,3\}\,$.}
\mp
To prove (5.8) consider the families $\,{\cal H}[\alpha]\;(\alpha<\omega)\,$
from above. Since $\,{\bf b}_{\omega_1}\,$ is the size of 
the family $\,{\cal U}:=
\bigcup_{\alpha<\omega_1}{\cal H}[\alpha]\,$, we obviously prove 
(5.8) by defining 
\sp
\qquad\qquad $\;{\cal H}'[\alpha]\;=\;
\{\,\Gamma_{5}\!\times\! H\;|\;H\in{\cal H}[\alpha]\,\}\,$\quad
for $\,\alpha<\omega_1\,$ and
\sp
\qquad\qquad $\;{\cal H}'[\omega_1]\;=\;
\{\,\Gamma_{7}\!\times\!\prod^*\!F\;\,|\,\;
F\in{\cal Y}({\cal U})\,\}\,.$
\mp\sp
Now we conclude the proof of Theorem 3 with the help of (5.8).
For $\,p\in\{2,3\}\,$ and 
each $\,\alpha\leq\omega_1\,$ define 
\mp
\cl{${\cal H}_p[\alpha]\,=\,
\{\,\Gamma_p\!\times\! H\;\,|\,\;H\in{\cal H}'[\alpha]\,\}\,.$}
\mp
Then both families 
$\,{\cal H}_2[\alpha]\,$ and $\,{\cal H}_3[\alpha]\,$
consist of precisely 
$\,{\bf b}_{\alpha+1}\,$ TFFI-groups of order $\,{\bf b}_{\alpha}\,$
for every $\,\alpha\leq\omega_1\,$. 
Furthermore, $\,{\cal H}_2[\alpha]\cap{\cal H}_3[\alpha']=\emptyset\,$
whenever $\,\alpha,\alpha'\leq\omega_1\,$. 
Since $\,\{2,3\}\cap\Omega(H)=\emptyset\,$ 
if $\,H\in{\cal H}'[\alpha]\,$ with $\,\alpha\leq\omega_1\,$,
the groups in 
$\,\bigcup_{\alpha\leq\omega_1}({\cal H}_2[\alpha]\cup{\cal H}_3[\alpha])\,$
are incomparable. 
\mp
For any ordinal $\,\alpha\,$ let
$\,\Lambda_\alpha\,$ denote the set of all cardinals $\,\kappa\,$
with $\,{\bf b}_\alpha<\kappa<{\bf b}_{\alpha+1}\,$.
Trivially, $\,|\Lambda_\alpha|\leq{\bf b}_{\alpha+1}\,$ for every $\,\alpha\,$. 
Under the General Continuum Hypothesis 
each set $\,\Lambda_\alpha\,$ would be empty.
However, by (7.1) in Section 7 we cannot 
rule out that $\,|\Lambda_0|={\bf b}_1\,$ and
$\,|\Lambda_{\alpha+1}|>{\bf b}_{\alpha+1}\,$ for every ordinal
$\,\alpha\,$. (Notice that $\,|\Lambda_0|={\bf b}_1\,$ is equivalent with 
$\;|\{\,|X|\;|\;X\subset\R\,\}|=|\R|\,$.)
\mp
Let $\,K\,$ be the set of all infinite cardinals smaller than 
$\,{\bf b}_{\omega_1}\,$ which are not beths. 
If $\,K=\emptyset\,$ then the proof of Theorem 3 is already finished 
by (5.8). So assume that $\,K\not=\emptyset\,$.
For each $\,\kappa\in K\,$ there is a unique ordinal 
$\,\alpha=\alpha(\kappa)<\omega_1\,$
such that $\,\kappa\in \Lambda_\alpha\,$.
Let $\,f\,$ be a function defined on $\,K\,$ 
such that $\,f(\kappa)\subset {\cal H}_3[\alpha(\kappa)]\,$
with $\,|f(\kappa)|=\kappa\,$ 
for every $\,\kappa\in K\,$
and $\,f(\kappa)\cap f(\lambda)=\emptyset\,$
whenever $\,\kappa,\lambda\in K\,$ are distinct. 
Such a function $\,f\,$ exists since $\,|\Lambda_\alpha|\,$ 
and $\,\kappa\in \Lambda_\alpha\,$ are not greater 
than $\,{\bf b}_{\alpha+1}=|{\cal H}_3[\alpha]|\,$ for every 
$\,\alpha<\omega_1\,$.
Finally put 
\mp
\qquad\qquad $\;{\cal R}[\kappa]\,=\,
\{\,\prod^*\! F\;\,|\,\;F\in{\cal Y}(f(\kappa))\,\}$\quad 
for every $\,\kappa\in K\,$ and  
\sp
\qquad\qquad $\;{\cal R}[{\bf b}_\alpha]\,=\,{\cal H}_2[\alpha]\;$\quad
for every $\,\alpha\leq \omega_1\,$
\mp
and define $\,{\cal H}\,$ as the union 
of the families $\,{\cal R}[\kappa]\,$ 
with $\,\kappa\,$ running through all infinite cardinals up to 
$\,{\bf b}_{\omega_1}\,$. 
A moment's reflection suffices to see 
that with this definition of $\,{\cal H}\,$ 
the proof of Theorem 3 is finished.
\bp
{\it Remark.} The upper bound $\,{\bf b}_{\omega_1}\,$ in Theorem 3
can be exceeded. If we work not only with the three sets 
$\,\P_1,\P_2,\P_3,\,$ but with infinitely many mutually disjoint infinite
sets $\,\P_i\subset\P\,$                        
then, in view of the previous considerations, 
the range $\,\aleph_0\leq\kappa\leq{\bf b}_{\omega_1}\,$ in Theorem 3 
can be expanded to $\,\aleph_0\leq\kappa\leq{\bf b}_{\omega_1\cdot\omega}\,$
(where $\,\omega_1\cdot\omega\,=\,\omega_1+\omega_1+\omega_1+\cdots\,$
with $\aleph_0$ summands).
Certainly, the new upper bound $\,{\bf b}_{\omega_1\cdot\omega}\,$ 
can also be exceeded. However, 
to cover the range $\,\aleph_0\leq\kappa\leq{\bf b}_{\omega_2}\,$
with the constructions of the previous proof 
is surely out of the question.
\bp\sp
{\bff 6. Proof of Theorem 2}
\mp
As any set of cardinal numbers, $\,K\,$ is well-ordered 
by the natural ordering.  While the {\it size} of $\,K\,$ 
is equal to the {\it cardinal number} $\,{\bf b}_{\omega_1}\,$, 
the {\it order type} of $\,K\,$ can be greater than 
the {\it ordinal number} $\,{\bf b}_{\omega_1}\,$.
(This is the reason why the inclusion of $\,\kappa={\bf b}_{\omega_1}\,$
in Theorem 3 is essential.)
Put $\,\vartheta={\bf b}_{\omega_1+1}\,$. Then 
$\,\vartheta\geq{\bf b}_{\omega_1}^+\,$ and hence
the order type of any subset of $\,K\,$ is less than 
$\,\vartheta\,$. Therefore we can define a set $\,\tilde K\,$ 
of {\it uncountable} cardinals such that 
$\,\tilde K\,$ contains $\,K\setminus\{\aleph_0\}\,$ 
and the order type of $\,\tilde K\,$ is $\,\vartheta\,$. 
Let $\,f\,$ be the unique strictly increasing 
mapping from the set of all ordinals $\,\alpha<\vartheta\,$
onto the well-ordered set $\,\tilde K\,$.
Naturally, $\,f(\alpha)\geq \alpha\,$ 
for every $\,\alpha<\vartheta\,$. 
\mp        
With $\,\Omega(\cdot)\,$ defined as in Section 5
let $\,{\bf G}\,$ be the class of all torsion-free 
groups $\,G\,$ where $\,\Omega(G)\not=\emptyset\,$
and $\,\{2,3\}\cap\Omega(G)=\emptyset\,$.
Let  $\,{\bf H}\,$ be the class of all 
groups which are not torsion-free 
and whose elements of infinite 
order are not $p$-divisible for any prime $\,p>3\,$.
Obviously, no group in the class $\,{\bf G}\,$ is embeddable in 
any group in the class $\,{\bf H}\,$ and vice versa.
\mp
In view of (4.1) and the remark in Section 4, 
for every $\,\kappa\in \tilde K\,$ 
there exists a family $\,{\cal H}[\kappa]\subset{\bf H}\,$ 
of $\,2^\kappa\,$ incomparable F-universal FI-groups of order $\,\kappa\,$. 
\mp
By (5.8) in the proof of Theorem 3, for 
every ordinal $\,\beta\leq \omega_1\,$ 
there exists a family $\,{\cal G}[{\bf b}_\beta]\subset{\bf G}\,$ 
of size $\,{\bf b}_{\beta+1}\,$ such that every 
$\,G\in{\cal G}[{\bf b}_\beta]\,$
is a FI-group of order $\,{\bf b}_\beta\,$
and the groups in 
$\,\bigcup\{\,{\cal G}[{\bf b}_\beta]\;|\;\beta\leq\omega_1\,\}\,$
are incomparable.
\mp
Since for infinite cardinals $\,\kappa,\lambda\,$ with $\,\kappa<\lambda\,$
there are precisely $\,\lambda\,$ ordinals between 
$\,\kappa\,$ and $\,\lambda\,$,
for each ordinal number $\,\beta\leq\omega_1\,$ 
we can choose a bijection $\,g_\beta\,$ from the ordinal interval
$\,\{\,\alpha\;|\;{\bf b}_\beta\leq \alpha <{\bf b}_{\beta+1}\,\}\,$
onto $\,{\cal G}[{\bf b}_\beta]\,$ and define 
\mp
\cl{${\cal F}[f(\alpha)]\,=\,\{\,g_\beta(\alpha)*H\;\,|\,\; 
H\in{\cal H}[f(\alpha)]\,\}$}
\mp
whenever $\,{\bf b}_\beta\leq\alpha<{\bf b}_{\beta+1}\,$ and $\,\beta\leq\omega_1\,$.
\mp
For $\,{\bf b}_\beta\leq\alpha<{\bf b}_{\beta+1}\,$ 
and $\,H\in{\cal H}[f(\alpha)]\,$
we have $\;|g_\beta(\alpha)*H|=\max\{|g_\beta(\alpha)|,|H|\}=
\max\{{\bf b}_\beta,f(\alpha)\}=f(\alpha)\;$
since $\,f(\alpha)\geq\alpha\,$.
Thus for every $\,\kappa\in \tilde K\,$
the family $\,{\cal F}[\kappa]\,$ consists of 
$\,2^\kappa\,$ groups of order $\,\kappa\,$. 
\mp
First we point out that the groups in $\,{\cal F}[\kappa]\,$ are incomparable
whenever $\,\kappa\in\tilde K\,$.
This is true in view of Corollary 2 
because if $\,G\,$ is a FI-group in $\,{\bf G}\,$ 
and $\,H_1,H_2\in{\cal H}[\kappa]\,$ are distinct
then the incomparable 
FI-groups $\,H_i\,$ are not embeddable in the FI-group $\,G\,$. 
\mp
Secondly, if $\,\kappa,\lambda\in \tilde K\,$ 
and $\,\kappa<\lambda\,$
then no group $\,F_\kappa\,$ (of order $\,\kappa\,$) in the family  
$\,{\cal F}[\kappa]\,$ is embeddable in a group $\,F_\lambda\,$ 
(of order $\,\lambda\,$) in the family  
$\,{\cal F}[\lambda]\,$ because the {\bf G}-factor of 
$\,F_\kappa\,$ is not embeddable in the {\bf G}-factor 
or in the {\bf H}-factor of $\,F_\lambda\,$.
\mp
Thus for every set $\,L\,$ 
with $\;K\subset L\subset\tilde K\;$
the definition  $\;{\cal F}\,=\,\bigcup\{\,{\cal F}[\kappa]\;|\;
\kappa\in L\,\}\;$ concludes the proof of Theorem 2
provided that $\,\aleph_0\not\in K\,$.
\mp
The case $\,\aleph_0\in K\,$ can be settled by adding the 
$\,2^{\aleph_0}\,$ incomparable groups 
$\,\Gamma[P]\times F\,$ from (3.2).
Indeed, these countable groups are F-universal FI-groups 
and it is evident that these groups are not embeddable in a {\bf G}-group 
or in a {\bf H}-group, whence they cannot be embedded 
in any group from the family 
$\;\bigcup\{\,{\cal F}[\kappa]\;|\;
\kappa\in \tilde K\,\}\,$.
This concludes the proof of Theorem 2.  
\bp\bp
{\bff 7. Long gaps between the beths}
\mp\sp
In the following, $\,0\in\N\,$.
Let $\,\kappa,\lambda,\mu\,$ always denote infinite
cardinal numbers and put $\;[\kappa,\lambda]\,=\,
\{\,\mu\;|\;\kappa\leq\mu\leq\lambda\,\}\,$.
It goes without saying that always $\,|[\kappa,\lambda]|\leq\lambda\,$.
Therefore, the following observation 
implies that {\it it is consistent 
with {\rm ZFC} that for every limit ordinal 
$\,\eta\,$ there exist $\,{\bf b}_\eta\,$ cardinal numbers 
up to $\,{\bf b}_\eta\,$.} 
(Of course, $\,\lambda\,$ is a strong limit 
if and only if $\,\lambda={\bf b}_\eta\,$ for a limit ordinal 
$\,\eta\,$.)
\mp
(7.1) $\;$ {\it It is consistent 
with {\rm ZFC} that 
the continuum function $\,\kappa\mapsto 2^\kappa\,$ is strictly increasing and
$\,|[{\bf b}_{n},{\bf b}_{n+1}]|={\bf b}_{n+1}\,$
for every even $\,n\in\N\,$
and $\,|[{\bf b}_{\beta+n},{\bf b}_{\beta+n+1}]|={\bf b}_{\beta+n+1}\,$
for every countable limit ordinal $\,\beta>0\,$
and every odd $\,n\in\N\,$ and 
$\,{\bf b}_{\alpha+1}<|[{\bf b}_{\alpha+1},{\bf b}_{\alpha+2}]|\,$ 
for every ordinal number $\,\alpha\,$.}
\mp
In order to proof (7.1) we apply a remarkable consistency theorem.
Let $\,\alpha\,$ be an ordinal number.
If $\,\kappa={\bf b}_\alpha^+\,$ 
then $\,\kappa\,$ is regular and, trivially,
$\,{\bf b}_{\alpha+1}\leq 2^{\kappa}\leq{\bf b}_{\alpha+2}\,$.
Therefore (and since $\,\aleph_0\,$ is regular and 
$\,{\bf b}_\beta\,$ is singular for every countable limit ordinal 
$\,\beta>0\,$) we obtain (7.1) as an immediate consequence of the following  
 fact.
\mp
(7.2)$\;$ {\it It is consistent with {\rm ZFC} that $\,\mu<\lambda\,$
implies $\,2^\mu<2^\lambda\,$ 
and $\;|[\aleph_0,2^\kappa]|=2^\kappa\;$  for every 
regular $\,\kappa\,$ 
and $\,2^\mu$ is regular for every 
singular $\mu$ and $\,2^\kappa$ is singular for every 
regular $\,\kappa\,$.}
\mp
Referring to Jech's profound text book [1],
a proof of (7.2) can be carried out as follows.
In G\"odel's universe L, let $\,\Theta(\kappa)\,$ denote 
the smallest fixed point of the $\aleph$-function 
whose cofinality equals $\,\kappa^+\,$.
Then in every generic extension of L which preserves 
cardinalities and cofinalities we have
$\;\forall\kappa:\,|[\aleph_0,\Theta(\kappa)]|=\Theta(\kappa)\;$ 
and hence $\;\forall\kappa:\,\Theta(\kappa)\not=\kappa^+\;$
due to $\;|[\aleph_0,\mu^+]|\leq\mu<\mu^+\;$ for all $\,\mu\,$.
\sp
By applying Easton's theorem [1, 15.18] one can 
create an Easton universe E generically extending L
such that the continuum function $\;\kappa\mapsto 2^\kappa=\kappa^+\;$
in L is changed into $\;\kappa\mapsto 2^{\kappa}=g(\kappa)\;$ in E
with $\,g(\kappa)=\Theta(\kappa)\,$ for every regular cardinal $\,\kappa\,$.
\sp
So in E we count $\,2^\kappa\,$ infinite cardinals below $\,2^\kappa\,$
for every regular $\,\kappa\,$. Furthermore, 
since $\,\Theta(\cdot)\,$ is strictly 
increasing in L, in E we have 
$\,2^\alpha<2^\beta\,$ whenever $\,\alpha,\beta\,$ are {\it regular} cardinals
with $\,\alpha<\beta\,$. Therefore and 
since the {\it Singular Cardinal Hypothesis} holds in E
(see [1] Ex.~15.12) and by [1] Theorem 5.22, 
 if $\,\mu\in{\rm E}\,$ is singular 
then $\,2^\mu\in{\rm E}\,$ is a successor cardinal and hence regular, while 
if $\,\kappa\in{\rm E}\,$ is regular then 
$\,2^\kappa\,$ equals $\,\Theta(\kappa)\,$ in E which 
is singular since $\,\kappa^+\,$ is the cofinality 
of $\,\Theta(\kappa)\,$ and $\,\Theta(\kappa)\not=\kappa^+\,$.
Consequently, the continuum function in E 
is strictly increasing on {\it all} cardinals.
This concludes the proof of (7.2).
\bp\bp\mp
{\bff References} 
\medskip\smallskip
[1] Jech, T.: {\it Set Theory.} 3rd ed. Springer 2002. 
\smallskip
[2] Kurosch, A.G.: {\it Gruppentheorie I.} Akademie-Verlag Berlin 1970.
\smallskip
[3] Macintyre, A., and Shelah, S.: {\it Uncountable
universal locally finite groups.}
J.~Alg.~{\bf 43}, 

\rightline{168-175 (1976).} 

[4] Rotman, J.: {\it An Introduction to the Theory of Groups.}  
Springer 1995. 
\smallskip
[5] Shelah, S.: {\it Infinite abelian groups, Whitehead problem and 
some constructions.}  

\rightline{Israel J.~Math.~{\bf 18}, 243-256 (1974).}

\bp\bp
{\sl Author's address:} Institute of Mathematics,
BOKU University, Vienna, Austria. 
\sp
{\sl E-mail:} {\tt gerald.kuba@boku.ac.at}
\bp\bp\bp
\hrule
\bp\bp\bp
{\bf This paper will be published 
in JOURNAL OF GROUP THEORY.}
\end